\documentstyle[11pt]{article}
\newtheorem{theorem}{Theorem}

\def\QED{\quad\blackslug\lower 8.5pt\null}

\begin{document}

\begin{center}
{\Large \bf  ON GEOMETRY OF HYPERSURFACES} 
\vspace*{2mm}

{\Large \bf  OF A PSEUDOCONFORMAL SPACE}
\vspace*{2mm}

{\Large \bf   OF LORENTZIAN SIGNATURE}

\vspace*{3mm}
{\large M.A. Akivis and  V.V. Goldberg}

\end{center}

\vspace*{5mm}

\begin{abstract}
There are three types of hypersurfaces in a pseudoconformal space 
$C^n_1$ of Lorentzian signature: spacelike, timelike, and 
lightlike. These three types of hypersurfaces are considered 
in  parallel. Spacelike hypersurfaces are endowed with a proper 
conformal structure, and timelike hypersurfaces 
are endowed with a conformal structure of Lorentzian type. 
Geometry of these two types of hypersurfaces can be studied 
in a manner that is similar to that for hypersurfaces of 
a proper conformal space. Lightlike hypersurfaces are endowed 
with a degenerate conformal structure. This is the reason that 
their investigation has special features. It is proved that under 
the Darboux mapping such hypersurfaces are transferred 
into tangentially degenerate $(n-1)$-dimensional submanifolds 
of rank $n-2$ located on the Darboux hyperquadric. 
The isotropic congruences of the space $C^n_1$ 
that are closely connected with lightlike hypersurfaces 
 and their Darboux mapping are also considered.
\end{abstract}
\setcounter{equation}{0}

\setcounter{section}{0}

{\bf 0. Introduction.}  
Submanifolds  in a proper conformal space $C^n$ were 
considered in numerous papers.  Submanifolds in pseudo-Euclidean 
spaces, in particular, in the Minkowski space, were also 
investigated in great detail (see, for example, [ON 83]). 
There are three types of submanifolds in a pseudo-Euclidean 
space:  spacelike, timelike, and lightlike. These 
three types of submanifolds were also studied in 
pseudo-Riemannian spaces of different signatures. In the recently 
published book [DB 96] the geometry of lightlike hypersurfaces 
and lightlike submanifolds of higher codimension in 
semi-Riemannian (or in another terminology pseudo-Riemannian) 
spaces were investigated in detail.

However, the property of submanifolds to be spacelike, timelike 
or lightlike is invariant with respect to conformal 
transformations of the pseudo-Riemannian spaces 
in which they are embedded. This is the reason that 
it is appropriate to consider all three types of submanifolds 
(spacelike, timelike, and lightlike) in the framework of 
pseudoconformal structures.

In the present paper we study hypersurfaces in a pseudoconformal 
space $C^n_1$ of Lorentzian signature. We show that the local 
theory of 
spacelike and timelike hypersurfaces in the space $C^n_1$ can be 
developed along the same lines as the theory of  hypersurfaces in 
a proper conformal space (Sections {\bf 4} and {\bf 5}). 
The theory of  lightlike 
(isotropic) hypersurfaces is quite different from  the theory of  
hypersurfaces in a proper conformal space. We consider some 
aspects of the theory of  lightlike  hypersurfaces (Section 
{\bf 6}) and  isotropic congruences  that are closely 
connected with  lightlike hypersurfaces (Section {\bf 7}). 
The use of pseudoconformal setting for studying of 
hypersurfaces allows us to apply the 
Darboux mapping, prove that under this mapping the image of 
a lightlike hypersurface is a tangentially degenerate submanifold 
in a projective space and describe singular points on a 
lightlike hypersurface and on an isotropic congruence of a 
pseudoconformal space.

Note that in  [DB 96] the results on 
lightlike hypersurfaces in semi-Riemannian spaces 
 are applied to electrodynamics and general relativity. 
But since many of these applications and 
the lightlike hypersurfaces themselves are conformally invariant, 
the results of the current paper can be used in similar and 
possibly other physical applications. 

The isotropic congruences in pseudo-Riemannian spaces 
are of interest for general relativity. In particular, 
they are connected with  construction of 
the Kerr metric describing black holes in the 
gravitational field (see [Ch 83], \$57).

{\bf 1. Preliminaries.}  
It is well known that a geometric model 
of spacetime in special relativity is the Minkowski space, 
that is, a four-dimensional pseudo-Euclidean space $R^4_1$ of 
signature $(3, 1)$ (see, for example, [BEE 96]). The fundamental 
quadratic form of this space is reduced to the form 
$$
g = (\omega^1)^2 + (\omega^2)^2 + (\omega^3)^2 - (\omega^4)^2,
$$
where  $\omega^1, \omega^2, \omega^3$ are the space coordinates, 
and $\omega^4$ is the time coordinate in the tangent space $T_x$ 
associated with a point $x$ of the space $R_1^4$. (The space $T_x$ is the set 
of vectors of the space $R_1^4$ emanating from 
the point $x$.) The relatively invariant fundamental form $g$ 
defines the Lorentzian 
metric in $R^4_1$. The group of transformations of the space 
$T_x$ preserving this metric is the pseudoorthogonal group 
${\bf O} (3, 1)$ that is also called the Lorentz group.

The equation $g = 0$ determines in the space $T_x$ the light 
cone 
$C_x$ whose generators are light trajectories propagating from a 
source located at the point $x$. The group of transformations of 
the space $T_x$ leaving the cone $C_x$ invariant is the group 
$G = {\bf O} (3, 1) \times {\bf H}$, where ${\bf H}$ is 
the group of homotheties of 
$T_x$.

Many results of special relativity, especially results concerning 
the light propagation, are connected with the conformal 
structure of the space $R^4_1$---the structure determined on 
$R^4_1$ by the relatively invariant 
fundamental form $g$. In fact, the equation 
$g = 0$ defines in $T_x$ the light cone $C_x$ with vertex at the 
point $x$, and the set of these cones is invariant under 
pseudoconformal transformations of the space $R^4_1$. Besides 
the light cones, these transformations leave invariant the 
set of hyperspheres of the space $R^4_1$ defined in $T_x$ by the 
equation $g = r^2$, where the number $r^2$ 
can be not only positive 
(as in the Euclidean space) but also negative or zero.
For $r^2 < 0$, the equation $g = r^2$ determines hyperspheres of 
imaginary radius; for $r^2 > 0$, 
it determines hyperspheres of real 
radius; and  for $r^2 = 0$, it determines hyperspheres of zero 
radius coinciding with the light cones \nopagebreak 
(see Figure 1 for $n = 3$). \pagebreak

\vspace*{3in} 
\begin{center}
Figure 1
\end{center} 
\vspace*{3mm}

Conformal transformations of the space $R^4_1$ form a group 
depending on 15 parameters. However, this group does not act 
bijectively in the space $R^4_1$. To make this group to act 
bijectively on the set of points of the space $R^4_1$, we 
should enlarge this set by ideal elements: a point at infinity 
$y = \infty$ and the light cone $C_y$ with vertex at the point 
$y$. The enlarged space $R^4_1$ is denoted by $C^4_1$ and is 
called the {\em pseudoconformal space} of signature $(3,1)$, 
$C^4_1 = R_1^4 \cup C_y$. After  this enlargement, the noncompact  
space $R^4_1$ becomes the compact space $C^4_1$. This is the 
 reason that such an operation is called the {\em 
compactification} of the Minkowski space.

In what follows we will consider not only four-dimensional space 
$C^4_1$ but also $n$-dimensional spaces $C^n_1$ of Lorentzian 
signature for $n \geq 3$. The fundamental form $g$ defining a conformal 
structure of this space can be reduced to the form
\begin{equation}\label{eq:1}
g = (\omega^1)^2 +  \ldots + (\omega^{n-1})^2 -  (\omega^n)^2.
\end{equation}
The space $C^n_1$ admits a one-to-one point mapping 
onto a hyperquadric $Q^n_1$ of a projective space $P^{n+1}$. 
The equation of $Q^n_1$ can be reduced to the form
\begin{equation}\label{eq:2}
 (x^1)^2 + \ldots  + (x^{n-1})^2 -  (x^n)^2 + (x^0)^2 
- (x^{n+1})^2 = 0.
\end{equation}
The projective coordinates $x^0, \ldots , x^{n+1}$ of 
points of the space $P^{n+1}$ are called {\em polyspherical 
coordinates} of the elements (points and hyperspheres) of 
the space $C^n_1$ (see [Kl 26] or [AG 96]).

The quadratic form on the left-hand side of 
equation (2) determines the {\em scalar product} of 
elements of the space $C^n_1$. As usual, we denote this scalar 
product by $(\hspace*{3mm}, \hspace*{3mm})$. The scalar square 
of a point of the space $C^n_1$ is equal to 0, and it is negative 
for spacelike hyperspheres and positive for timelike 
hyperspheres. The vanishing of the scalar product 
of two hyperspheres means that the hyperspheres are orthogonal, and the 
vanishing of the scalar product of a point and a 
hypersphere means that the point belongs to the hypersphere.

The group of conformal transformations 
of the space $C^n_1$ is isomorphic to 
the group of projective transformations of the space 
$P^{n+1}$ sending the hyperquadric $Q^n_1$  to itself. 
This group is denoted   by ${\bf PO} (n+2, 2)$ 
and is expressed as follows:
$$
{\bf PO} (n+2, 2) := \left\{ 
\begin{array}{ll}
{\bf SO} (n+2, 2) & 
\mbox{if} \;\; n \;\; \mbox{is odd}, \\
 {\bf O} (n+2, 2)/{\bf Z}_2 & \mbox{if} \;\; n \;\; 
\mbox{is even},
\end{array}
\right.
$$
where ${\bf O} (n+2, 2)$ and ${\bf SO} (n+2, 2)$ are 
the groups of pseudoorthogonal and special pseudoorthogonal 
transformations of the indicated signature, respectively, 
and ${\bf Z}_2$ is the cyclic group of second order. 
In both cases this group depends on $\frac{1}{2} (n+1)(n+2)$ 
parameters.

The mapping $\varphi: C_1^n \rightarrow Q_1^n$ is 
called the {\em Darboux mapping}, and the hyperquadric $Q^n_1$ 
is called the {\em Darboux hyperquadric}. Such a mapping was 
constructed first for the proper conformal three-dimensional 
space $C^3$ (see [Kl 26]). Under the mapping $\varphi$ 
to the isotropic cones $C_x$ there correspond the {\em asymptotic 
cones} of the hyperquadric $Q^n_1$. This  hyperquadric 
carries real rectilinear generators to which in the space $C^n_1$ 
there correspond the lines of light propagation. The light 
cones in $C^n_1$ are called also the {\em isotropic cones}, and 
the lines of light propagation are called the {\em isotropic 
lines} of the space $C^n_1$.

Further we will apply the  method of moving frames. In the space 
$C^n_1$ we consider a family of  conformal frames consisting 
of two points  $A_0$ and $A_{n+1}$ and $n$ hyperspheres $A_r, 
r = 1, \ldots , n$, 
passing through these points. The frame elements of such frames 
satisfy the following analytical conditions: 
\begin{equation}\label{eq:3}
 (A_0, A_0) = (A_{n+1}, A_{n+1}) = 0, 
\end{equation}
\begin{equation}\label{eq:4}
\renewcommand{\arraystretch}{1.3}
\begin{array}{ll}
 (A_0, A_0) = (A_{n+1}, A_{n+1}) = 0, \;\;
  (A_0, A_r) = (A_{n+1}, A_r) = 0, \\
(A_r, A_s ) = g_{rs}, \;\;\;\;\; r, s = 1, \ldots , n.
\end{array}
\renewcommand{\arraystretch}{1}
\end{equation}
In addition, we normalize the points $A_0$ and $A_{n+1}$ by the 
condition 
\begin{equation}\label{eq:5}
 (A_0, A_{n+1}) = - 1.
\end{equation}

Under the Darboux mapping to such frames there correspond point 
projective frames in the space $P^{n+1}$ for which the points 
$A_0$ and $A_{n+1}$ lie on the Darboux hyperquadric but do not 
belong to any of its rectilinear generators, and the points $A_r$ 
form a basis of the $(n-1)$-dimensional subspace that is 
polar-conjugate to the straight line $A_0 A_{n+1}$ with respect 
to the Darboux hyperquadric. With respect to this projective 
point frame the equation of Darboux hyperquadric takes the form
\begin{equation}\label{eq:6}
g_{rs} x^r x^s - 2 x^0 x^{n+1} = 0,
\end{equation}
where the quadratic form $g_{rs} x^r x^s$ is of signature 
$(n-1, 1)$.

The equations of infinitesimal displacement of our conformal  
frame in the space $C^n_1$ are
\begin{equation}\label{eq:7}
dA_\xi = \omega_\xi^\eta A_\eta, \;\;\;\;\; \xi, \eta = 0, 1, \ldots, n + 1,
\end{equation}
where  $ \omega_\xi^\eta $ are differential 1-forms containing 
the parameters, on which the group ${\bf PO} (n+2, 2)$ 
depends, and their differentials.

If we differentiate conditions (3)--(5) by means of 
equations (7), we obtain that the forms $\omega_\xi^\eta $ 
satisfy the following equations:
\begin{equation}\label{eq:8}
\omega_{n+1}^0 = \omega_0^{n+1} = 0, \;\; 
\omega_0^0 + \omega_{n+1}^{n+1} = 0, 
\end{equation}
\begin{equation}\label{eq:9}
\omega_r^{n+1} - g_{rs} \omega_0^s = 0, \;\;  
\omega_r^0 - g_{rs} \omega_{n+1}^s = 0,
\end{equation}
\begin{equation}\label{eq:10}
dg_{rs} = g_{rt} \omega_s^t + g_{ts} \omega_r^t.
\end{equation}
In addition, the forms $\omega_\eta^\xi$ satisfy 
the structure equations of the spaces $C_1^n$ and $P^{n+1}$:
\begin{equation}\label{eq:11}
d \omega_\xi^\eta =  \omega_\xi^\zeta \wedge \omega_\zeta^\eta 
\end{equation}
that are necessary and sufficient conditions for complete 
integrability of equations (7).

{\bf 3. Hypersurfaces in the space 
\protect\boldmath \( C^n_1. \) \protect\unboldmath} In the space 
$C_1^n$ 
we consider a hypersurface $V^{n-1}$, that is, a smooth, 
connected and simply connected submanifold of dimension $n-1$. 
The conformal structure of the space $C_1^n$ induces a conformal 
structure on the hypersurface $V^{n-1}$. The nature of this 
structure depends on the mutual location of tangent 
hyperplanes $T_x (V^{n-1}) = \tau_x, x \in V^{n-1}$, with respect 
to the isotropic cones $C_x$ of the space $C^n_1$. Three 
``pure'' cases of such location are possible;

\begin{description}
\item[a)] At any point $x \in V^{n-1}$ the hyperplane $\tau_x$ 
and the isotropic cone $C_x$ have only one common point $x$. 
Then the quadratic form $\widetilde{g} = g|_{\tau_x}$ 
on the  hypersurface $V^{n-1}$ is positive definite, and 
on   $V^{n-1}$ a proper conformal structure is induced. 
A hypersurface $V^{n-1}$ of this type is called {\em spacelike}.

\item[b)] At any point $x \in V^{n-1}$ the hyperplane $\tau_x$ 
intersects  the isotropic cone $C_x$ along a real cone 
$\widetilde{C}_x$ of dimension $n-2$. Then the quadratic form 
$\widetilde{g}$ on  $V^{n-1}$ has  signature $(n-2, 1)$, and 
on   $V^{n-1}$ a conformal structure of the same signature 
is induced. A hypersurface $V^{n-1}$ of this type is called 
{\em timelike}.

\item[c)] At any point $x \in V^{n-1}$ the hyperplane $\tau_x$ 
is tangent to the isotropic cone $C_x$. Then the quadratic form 
$\widetilde{g}$ on  $V^{n-1}$ has signature $(n-2, 0)$. 
A hypersurface $V^{n-1}$ of this type is called {\em lightlike} 
or {\em isotropic}.
\end{description}

For the dimension $n = 3$ these three cases are represented 
on Figures 2, 3, and 4.

The terminology (spacelike, timelike and lightlike) 
is related to that of general relativity. 
As was noted in Introduction, 
spacetime in special  relativity is a Minkowski space. 
In general relativity it is a pseudo-Riemannian space. 
In both 

\begin{minipage}{2.5in}
\vspace{2.5in}
\begin{center}
Figure 2
\end{center} 
\end{minipage} 
\begin{minipage}{3.3in}
\vspace{2.5in} 
\begin{center}
Figure 3
\end{center} 
\end{minipage} 

\vspace*{5mm}
\vspace{2.5in} 
\begin{center}
Figure 4
\end{center} 
\vspace*{3mm}

\newpage

\vspace*{2.5in} 
\begin{center}
Figure 5
\end{center} 
\vspace*{3mm}

\noindent 
cases its metric has the signature $(3, 1)$ 
(or $(1, 3)$---this depends on the method of presentation). 
In general relativity the  isotropic cone $C_x$ plays the 
role of the light cone. This cone divides the tangent 
space $T_x (C^n_1)$ (or space $T_x (C^n_{n-1})$) into 
two domains---internal and external. Directions belonging 
to the first domain are called timelike, and directions belonging 
to the second domain are called spacelike (see Figure 5).
The tangent hyperplane $T_x (V^{n-1})$ to a spacelike 
hypersurface contains only directions located outside of the cone 
$C_x$, namely spacelike directions. For a timelike hypersurface 
$V^{n-1}$ the tangent hyperplane $T_x (V^{n-1})$ contains both 
spacelike and timelike  directions. 

Note that hyperspheres of real radius, defined in the space $T_x$ 
by the equation $g = a$ where $ a > 0$, are spacelike 
hypersurfaces without singularities. If $a < 0$, then the 
equation 
$g = a$ defines timelike hypersurfaces also not having 
singularities. Finally, if $a = 0$, then the equation 
$g = a$ defines a hypersphere of zero radius, that is, 
a lightlike hypersurface with the only singular point $x$. 
For $n = 3$, such hypersurfaces are presented in Figure 1.

Note also that although under conformal transformations 
hyperspheres are transferred into hyperspheres, 
the radii of these hyperspheres are not invariant. 
However, under conformal transformations the nature 
of hyperspheres (i.e., their property to be spacelike, or 
timelike or lightlike) is invariant. 

Besides ``pure'' hypersurfaces indicated above, there are 
hypersurfaces having points of two or of 
all three types indicated 
above. However, we will not consider such hypersurfaces 
in the present paper.

{\bf 4. Geometry of spacelike hypersurfaces.} 
We will study now the geometry of spacelike hypersurfaces 
$V^{n-1}$ of the pseudoconformal space $C^n_1$  in more detail. 

 With each point $x$ of the hypersurface $V^{n-1}$, we  associate  
a family of conformal frames in such a way that $A_0 = x$, the 
hypersphere $A_n$ is tangent to  $V^{n-1}$ at the point $x$, and 
the hyperspheres $A_i, \; i = 1, \ldots , n - 1$, 
 are orthogonal to  $V^{n-1}$ at this point. Hence,
 the hypersphere 
$A_n$ is spacelike, and hyperspheres $A_i$ are timelike.

After such a specialization of moving frames 
equations (3) and (5) will not be changed as well as the first 
two groups of equations (4) while the third group of 
equations (4) takes the form
\begin{equation}\label{eq:12}
 (A_i, A_n) = 0, \;\; (A_i, A_j ) = g_{ij}, \;\; 
(A_n, A_n) = - 1,
\end{equation}
where $(g_{ij})$ is a nondegenerate symmetric matrix of 
coefficients of a positive definite quadratic form $\widetilde{g}$. Note that 
the last equation in (12) 
is obtained by means of an additional normalization of 
the hypersphere $A_n$: this normalization is possible, 
since $A_n$ is a spacelike hypersphere.
 With respect to this frame the equation 
of the Darboux hyperquadric takes the form 
\begin{equation}\label{eq:13}
g_{ij} x^i x^j - (x^n)^2 - 2 x^0 x^{n+1} = 0.
\end{equation}

Since the hypersphere $A_n$ is tangent to the hypersurface 
$V^{n-1}$ at the point $A_0$, we have $(d A_0, A_n) = 0$. By the 
first equation of (7), this implies 
\begin{equation}\label{eq:14}
\omega^n_0 = 0    
\end{equation}
and
\begin{equation}\label{eq:15}
 d A_0 = \omega^0_0 A_0 + \omega^i A_i,
\end{equation}
where $\omega^i = \omega_0^i$. 
>From (15) it follows that the forms $ \omega^i$ are linearly 
independent. 

The family of frames described above is the bundle  
${\cal R}^1 (V^{n-1})$ of first-order 
frames associated with the hypersurface $V^{n-1}$. 
A base of this frame bundle is the hypersurface $V^{n-1}$, 
and its  fiber is the collection of frames with a fixed point 
$x = A_0$. The forms $\omega^i$ are {\em base forms} of   
${\cal R}^1 (V^{n-1})$, and the forms $\omega^0_0, \omega^0_i, 
\omega^j_i$ and $\omega^0_n$ are its {\em fiber forms}.

The quadratic form $g$, defining the conformal structure  
in the space $C_1^n$ at the point $x$, is expressed now as 
$$
g = g_{ij} \omega^i \omega^j - (\omega^n)^2,
$$
and its restriction to the hypersurface $V^{n-1}$ becomes 
\begin{equation}\label{eq:16}
\widetilde{g} = g_{ij} \omega^i \omega^j. 
\end{equation}
The form $\widetilde{g}$ is positive definite and defines a 
proper conformal structure on the hypersurface $V^{n-1}$.  
The coefficients $g_{ij}$ of this quadratic form 
generate a $(0, 2)$-tensor. This tensor is 
associated with a first-order neighborhood of 
the  hypersurface $V^{n-1}$, since by (15) we have
$$
\widetilde{g}= (dA_0, dA_0).
$$
We will not write here all equations which 
the forms $\omega^\eta_\xi$ satisfy 
on the  hypersurface $V^{n-1}$ and equations obtained 
as differential prolongations of equations (14). 
They differ unessentially from  similar equations 
in the theory of  hypersurfaces in a proper conformal 
space $C^n$. The latter theory was considered in 
[A 52] and [SS 80] (see also [AG 96]). 
 This is the reason that we are not going 
to consider in detail this construction as well as other topics 
of the theory of spacelike hypersurfaces which are known 
 for hypersurfaces of the proper conformal space $C^n$. 

{\bf 5. Geometry of timelike hypersurfaces.} 
Suppose now that a hypersurface $V^{n-1} \subset C^n_1$ is 
timelike. Then at any point $x \in V^{n-1}$ 
its tangent hyperplane $T_x (V^{n-1}) = \tau_x$ is located 
with respect to the  light cone $C_x$ as indicated in Figure 
3. The tangent hyperspheres to the hypersurface $V^{n-1}$ 
are timelike. Thus, they can be normalized by the condition
$$
(A_n, A_n) = 1.
$$
A timelike hypersurface $V^{n-1}$ is also determined by equation 
(14).

The fundamental form $g$ defining the conformal structure of 
the space $C^n_1$ is expressed now as
$$
g = g_{ij} \omega^i \omega^i + (\omega^n)^2,
$$
and its restriction  $\widetilde{g}$ to $V^{n-1}$ becomes 
$$
\widetilde{g} = g_{ij} \omega^i \omega^j. 
$$
However unlike for spacelike hypersurfaces, for the timelike 
hypersurfaces 
the form $\widetilde{g}$ is of signature $(n-2, 1)$. Thus,  
 {\em a timelike hypersurface $V^{n-1}$ possesses 
a pseudoconformal structure of Lorentzian signature}.

However, again  the system of equations associated with a 
timelike hypersurface $V^{n-1}$ differs  unessentially from  
similar equations in the theory of  hypersurfaces in a proper 
conformal space $C^n$. Thus, we will not go into details of 
investigation of timelike hypersurfaces.

Note only that since the  isotropic cone $\widetilde{C}_x$ of 
a timelike hypersurface is real, its mutual location with 
the cone $a_{ij} \omega^i \omega^j = 0$, determined by the 
second fundamental tensor $a_{ij}$ of $V^{n-1}$ 
and connected with a second-order neighborhood of 
a point $x \in V^{n-1}$, can be more diverse than for 
a hypersurface of the space $C^n$ or for a spacelike 
hypersurface of the space $C^n_1$. It would be interesting 
to construct a classification of timelike hypersurfaces 
based on the location of these two cones.

{\bf 6. Geometry of  lightlike hypersurfaces.}  Next we 
consider  lightlike hypersurfaces of the 
space $C^n_1$. For such hypersurfaces the quadratic form 
$\widetilde{g}$ is of signature $(n-2, 0)$, and 
they carry  degenerate conformal structures. 

Our considerations will be simpler if we consider the Darboux 
mapping of a  lightlike hypersurface $V^{n-1} \subset C_1^n$ and 
all geometric objects associated with this hypersurface. The  
hypersurface $V^{n-1} $ will be mapped onto a submanifold 
$U^{n-1}$ 

\newpage

\vspace*{2.0in} 
\begin{center}
Figure 6
\end{center} 
\vspace*{3mm}

\noindent 
of dimension $n - 1$ belonging to the Darboux 
hyperquadric that is determined in the space $P^{n+1}$ by the 
equation (6).

As usual we locate the vertex $A_0$ of the moving frame at the 
varying point $x \in U^{n-1}$ and the vertices $A_1, \ldots , 
A_{n-1}$ in the tangent $(n-1)$-plane $T_x (U^{n-1})$. Then 
the  equations 
\begin{equation}\label{eq:17}
\omega_0^n  = 0
\end{equation}
holds.

But since the hypersurface $V^{n-1}$ is  lightlike, the tangent 
$(n-1)$-plane $T_x (U^{n-1})$ is tangent to the asymptotic 
cone $C_x$ of the Darboux hyperquadric. The latter cone 
corresponds to the  isotropic cone $C_x$ 
of the space $C^n_1$. We place the vertex 
$A_1$ on the rectilinear generator along which the cone $C_x$ is 
tangent to the subspace $T_x (U^{n-1})$. We also place the vertex 
$A_n$ on the cone $C_x$ but outside of this tangent subspace 
$T_x (U^{n-1})$ (see Figure 6). 
Then in addition to equations (3)--(5) which the elements of a 
moving frame satisfy, we have also the following relations: 
\begin{equation}\label{eq:18}
 (A_1, A_1) = (A_n, A_n) =  (A_0, A_1) =  (A_0, A_n) = 0. 
\end{equation}
Moreover, we normalize the points $A_1$ and $A_n$ by the 
condition
\begin{equation}\label{eq:19}
  (A_1, A_n) = -1.
\end{equation}

By virtue of this, the matrix of the 
scalar products of the elements 
of the moving frame takes the form 
\begin{equation}\label{eq:20}
  (A_\xi, A_\eta) = \pmatrix{ 
0 & 0 & 0 & 0 & -1 \cr
0 & 0 & 0 & -1 & 0 \cr
0 & 0 & g_{ij} & 0 & 0 \cr
0 & -1 & 0 & 0 & 0 \cr
-1 & 0 & 0 & 0 & 0 \cr},
\end{equation}

\newpage 

\vspace*{1.7in} 
\begin{center}
Figure 7
\end{center} 
\vspace*{3mm}

\noindent 
where $\xi, \eta = 0, 1, \ldots , n+1; i, j = 2, \ldots , n-1$. 
As a result, the equation of the Darboux hyperquadric takes 
the form
\begin{equation}\label{eq:21}
g_{ij} x^i x^j  - 2 x^1 x^n - 2 x^0 x^{n+1} = 0,
\end{equation}
where $g_{ij} x^i x^j $ is a positive definite quadratic form. 

It follows that the $(n-3)$-dimensional subspace, 
determined in the space $P^{n+1}$ by the points $A_i, \; i = 2, 
\ldots , n - 1$, does not have real common points with the 
Darboux hyperquadric, and the subspace, which is polar-conjugate 
to the above subspace with respect to this hyperquadric and 
is determined by the points $A_0, A_1, A_n$ and $A_{n+1}$, 
intersects this hyperquadric in the following ruled surface of 
second order:
$$
x^1 x^n + x^0 x^{n+1} = 0.
$$
The above four points are located on this ruled 
surface as indicated in Figure 7.

The equation of the asymptotic cone $C_x$ at the point 
$x = A_0$ of the Darboux hyperquadric has the form 
\begin{equation}\label{eq:22}
g = g_{ij} \omega^i \omega^j  - 2 \omega^1 \omega^n = 0.
\end{equation}
Since the equation of 
the hypersurface $V^{n-1}$ has the form (17), 
the equation of the   cone $\widetilde{C}_x$ of the 
submanifold $U^{n-1}$ as well as of that of  the 
hypersurface $V^{n-1}$ has the form 
\begin{equation}\label{eq:23}
\widetilde{g} = g_{ij} \omega^i \omega^j  = 0,\;\; i, j = 2, 
\ldots , n-1.
\end{equation}
This implies that  at the point $x$ this  light  cone has 
a single rectilinear generator $A_0 A_1$ along which the subspace 
$T_x (U^{n-1})$ is tangent to the  asymptotic cone $C_x$. 

Next we write the equations of infinitesimal 
displacement of the moving frame associated with the point 
$x \in U^{n-1} \subset Q^n_1 \subset  P^{n+1}$ in the form (7) 
where the 1-forms $ \omega_\xi^\eta$ satisfy  
the equations (8)--(10)and also 
the equations  obtained by differentiation of equations 
(18) and (19): 
\begin{equation}\label{eq:24}
\omega_1^n= 0, \;\; \omega_n^1= 0, \;\; \omega_0^n 
+ \omega_1^{n+1}= 0, \;\; \omega_0^1 + \omega_n^{n+1}= 0, \;\;
\omega_1^1 + \omega_n^n= 0.
\end{equation}
If we also differentiate the equation $g_{1i} = 0$, we find that
$$
\omega_i^n= g_{ij} \omega_1^j, \;\; i, j = 2, 
\ldots , n-1.
$$
Since the tensor $g_{ij}$ is nondegenerate, it follows from 
the last equation that 
\begin{equation}\label{eq:25}
\omega_1^i = g^{ij}  \omega_j^n.
\end{equation}

Next, taking exterior derivatives of equation (17) and 
taking into account the first equation of (24), we obtain
\begin{equation}\label{eq:26}
\omega_i^n \wedge   \omega_0^i = 0, \;\;\;\;\; i = 2, \ldots , 
n -1.
\end{equation}
 Applying Cartan's lemma to equation (26), we find 
that 
$$
\omega_i^n = \lambda_{ij}   \omega_0^j, \;\; i, j = 2, 
\ldots , n-1,
$$
where $\lambda_{ij} = \lambda_{ji}$. Taking into account 
equations (25), we find 
\begin{equation}\label{eq:27}
\omega_1^i = g^{ik} \lambda_{kj}   \omega_0^j 
= \lambda_j^i \omega_0^j,
\end{equation}
where $\lambda_j^i = g^{ik} \lambda_{kj}$ is a symmetric 
nondegenerate affinor.

We consider now the differentials of the points $A_0$ and 
$A_1$. By (17) and (7), we obtain 
\begin{equation}\label{eq:28}
\left\{
\begin{array}{ll}
dA_0 = \omega_0^0 A_0 +  \omega_0^1 A_1 + \omega_0^i A_i, \\
dA_1 = \omega_1^0 A_0 +  \omega_1^1 A_1 + \omega_1^i A_i.
\end{array}
\right. 
\end{equation}
 From  equations (27) and (28) it follows that if $\omega_0^i 
= 0$, then the point $A_0$ moves along the  lightlike straight 
line $A_0 A_1$ belonging to the cone $C_x$ and  describes the 
entire line $A_0 A_1$. This means that the 
submanifold $U^{n-1}$ is a ruled submanifold. Moreover, the 
1-form $\omega_0^1$ defines the displacement of the point $A_0$ 
along the straight line $A_0 A_1$.

Next, equations (28) show that at any point of the straight 
line $A_0 A_1$,  the tangent $(n-1)$-dimensional subspace  is 
fixed and coincides with the subspace 
$T_x (U^{n-1}) = A_0 \wedge A_1 \wedge A_2  \wedge \ldots  \wedge  
A_{n-1}$. Thus,   the submanifold $U^{n-1}$ is tangentially 
degenerate of rank $n - 2$ (see [AG 93], Ch. 4), 
since the tangent subspace 
$T_x (U^{n-1})$ depends precisely on $n - 2$ parameters. 

Let $X = A_1 + x A_0$ be an arbitrary point of the rectilinear 
generator $A_0 A_1$ of the submanifold $U^{n-1}$. Its 
differential is determined by the formula
$$
dX \equiv (\omega_1^i + x \omega_0^i) A_i \pmod{A_0, A_1}.
$$
Since, by (27), 
$$
\omega_1^i + x \omega_0^i = (\lambda_j^i + x \delta_j^i) 
\omega_0^j,
$$
there are singular points on the straight line $A_0 A_1$, 
and their coordinates are determined by the equation 
\begin{equation}\label{eq:29}
\det (\lambda_j^i + x \delta_j^i) = 0.
\end{equation}
Since the tensor $\lambda_j^i$ is symmetric, this equation 
has $n - 2$ real roots if we count each root as many times as its 
multiplicity.
 
Thus, we have proved the following result.

\begin{theorem}
Under the Darboux mapping, to a lightlike hypersurface $V^{n-1}$ 
of the pseudoconformal space $C^n_1$ there corresponds a ruled 
tangentially degenerate submanifold $U^{n-1}$ of rank $n-2$ whose 
rectilinear generator carries $n-2$ real singular points if each 
of them is counted as many times as its multiplicity. These 
points are the images of singular points of the  lightlike 
hypersurface $V^{n-1}$. 
\end{theorem}

The loci of singular points on  lightlike hypersurfaces $V^{n-1}$  
are submanifolds whose dimension is less than $n - 1$. These 
submanifolds are called {\em focal submanifolds}. The dimension 
of  focal submanifolds depends on the multiplicity of  
their elements---singular points.

If $x_1$ is a simple root of equation (29), then  to this 
root there corresponds a family of torses (developable surfaces) on the 
submanifold $U^{n-1}$ which are defined by the system of 
equations 
\begin{equation}\label{eq:30}
\omega_1^i + x_1 \omega_0^i = 0.
\end{equation}
>From the well-known theorem of linear algebra on orthogonality of 
eigendirections of a symmetric linear operator, it follows that 
to distinct roots of equation (29) there correspond two 
mutual orthogonal families of torses on $U^{n-1}$. It is not 
difficult also to describe  submanifolds on $U^{n-1}$ 
corresponding to multiple roots of equation (29).

Note that since in general relativity, to  lightlike straight 
lines of the space $C_1^4$ there correspond lines of propagation 
of light, then to singular points on  lightlike 
hypersurfaces there correspond sources of light or points of its 
absorption, and their focal submanifolds are lighting surfaces or 
surfaces of light absorption. The further study of  lightlike 
hypersurfaces in the space $C_1^4$ can be of interest for general 
relativity. 

Note that the theory of lightlike hypersurfaces in 
semi-Riemannian spaces was studied in detail in [DB 96] 
and that some problems of the global theory of such hypersurfaces 
were considered by Kossowski (see, for example, [K 89]). 

{\bf 7.  Isotropic congruences.}  
The notion of  isotropic congruences of the space $C^n_1$ is 
closely connected with the theory of  isotropic hypersurfaces. 
An {\em  isotropic congruence} is an $(n-1)$-parameter family of 
isotropic straight lines such that through a generic point lying 
in a sufficiently small neighborhood of a straight line of the 
family there passes a unique straight line of the family. 

To study the isotropic congruences we will apply again 
the Darboux mapping of the space $C^n_1$. 
Consider the set $U$ of rectilinear generators of 
the Darboux hyperquadric $Q_1^n$. With any rectilinear 
generator of $Q_1^n$ we associate a family of frames 
described in Section {\bf 6}. Then 
with respect to any such frame the Darboux hyperquadric 
is defined by equation (21), and 
the components of infinitesimal displacements of these frames 
satisfy  equations (8)--(10) and (24).

Consider the rectilinear generator  $A_0 A_1$. We have 
\begin{equation}\label{eq:31}
dA_0 = \omega_0^0 A_0 + \omega_0^1 A_1 + \omega_0^i A_i + \omega_0^n A_n 
\end{equation}
and
\begin{equation}\label{eq:32}
dA_1 = \omega_1^0 A_0 + \omega_1^1 A_1 + \omega_1^i A_i - \omega_0^n A_{n+1}, 
\end{equation}
where $i = 2, \ldots , n - 1$. 
On the hyperquadric $Q_1^n$ the forms $\omega_0^i, \omega_1^i$, 
and $\omega_0^n$ are linearly independent, and their number is 
equal to $2n - 3$. Thus {\em the set $U$ is a differentiable 
manifold of dimension $2n-3$}.

The congruence of isotropic straight lines is 
an $(n-1)$-dimensional submanifold $S$ of the manifold $U$. 
In general, this submanifold can be given on $U$ by the 
following system of $n - 2$ equations:
\begin{equation}\label{eq:33}
\omega^i_1 = \lambda_j^i \omega_0^j + \lambda^i \omega^n_0, 
\;\;\;\; i, j = 1, 2, \ldots , n - 1.
\end{equation}
The forms $\omega_0^i$ and $\omega_0^n$ are 
basis forms of the congruence $S$.

On the congruence $S$ equation (32) takes the form
\begin{equation}\label{eq:34}
dA_1 = \omega_1^0 A_0 + \omega_1^1 A_1 
+ \lambda_j^i \omega^j_0 A_i -  \omega_0^n (A_{n+1} 
- \lambda^i A_i)
\end{equation}

In the projective space $P^{n+1}$ the straight lines of 
the congruence in question describe a hypersurface that we will 
also denote by $S$. As a point set, the hypersurface 
$S$ coincides with an open domain of the hyperquadric $Q_1^n$. 

Let us study properties of the  hypersurface $S$. 
Equations (31) and (34) imply that the linear span of a 
first-order neighborhood of 
the generator $A_0 A_1$ coincides with the entire space 
$P^{n+1}$. Next all tangent hyperplanes $T_x (S)$ at the points 
of its generator $A_0 A_1$ have  the common subspace 
$A_0 \wedge A_1 \wedge A_2 \wedge\ldots \wedge A_{n-1}$.

Consider singular points of the  hypersurface $S$. 
Its point $X = A_1 + x A_0$ is  singular if at this point the 
dimension of the tangent subspace $T_x (S)$ is less than $n$. 
By (31) and (34) we have 
$$
dX \equiv 
(\lambda_j^i + x \delta_j^i) \omega_0^j A_i 
+ \omega_0^n (x A_n - A_{n+1} + \lambda^i A_i) \pmod{A_0, A_1}.
$$
Thus, the dimension 
of the tangent subspace $T_x (S)$ is less than $n$ if and only if 
\begin{equation}\label{eq:35}
\det (\lambda_j^i + x \delta_j^i) = 0.
\end{equation}
Equation (35) determines singular points on the hypersurface $S$. 
Equation (35) differs from equation (29), determining singular 
points on  a lightlike  hypersurface of the space $C_1^n$, 
only by 
the fact that the affinor $\lambda_j^i$ was symmetric in (29) and 
is not symmetric in (35). 

Now suppose that the equation $\omega^n_0 = 0$ is 
completely integrable. Then the hypersurface $S$ is stratifited 
into a one-parameter family of $(n-1)$-dimensional 
submanifolds to which in the space $C^n_1$ there correspond lightlike 
hypersurfaces $V^{n-1}$. 

The condition of complete integrability of the equation 
$\omega^n_0 = 0$ has the form  
$d \omega^n_0 \wedge \omega_0^n = 0$. 
By structure equations (11),  
this implies that in equation (33)  
the affinor $\lambda_j^i$ is symmetric. As a result, all 
singular points of a rectilinear generator of the ruled 
hypersurface $S$ are real. 

 Thus the following theorem is valid:

\begin{theorem}
Any rectilinear 
generator $A_0 A_1$ of the isotropic congruence $S$ 
carries $n-2$  singular points if each 
of them is counted as many times as its multiplicity, and 
some of these singular points nor all of them can be complex. 
These singular points are the Darboux images of the 
singular points of the congruence of the space $C^n_1$. 
 If equation $\omega_0^n = 0$ is completely integrable on $S$, 
then it determines a stratification of 
the congruence $S$ into lightlike hypersurfaces,  and all these 
singular points are real.
\end{theorem}

{\em Authors' addresses}:\\

\begin{tabular}{ll}
M.A. Akivis &                           V.V. Goldberg \\
Department of Mathematics      &      Department of Mathematics\\
Ben-Gurion University of the Negev & 
                       New Jersey Institute of Technology \\
P.O. Box 653 & University Heights \\
Beer Sheva 84105, Israel & Newark, NJ 07102, U.S.A.     
\end{tabular}

\end{document}